\documentclass[11pt]{article}

\usepackage{latexsym}
\usepackage{amssymb}

\newtheorem{theorem}{Theorem}

\newtheorem{definition}{Definition}

\newtheorem{example}{Example}

\begin{document}
{
\begin{center}
{\Large\bf
On modified kernel polynomials and classical type Sobolev orthogonal polynomials}
\end{center}
\begin{center}
{\bf S.M. Zagorodnyuk}
\end{center}

\section{Introduction.}

The theory of orthogonal polynomials on the real line has a great amount of contributions and a lot of applications in various areas of science, 
see~\cite{cit_50000_Gabor_Szego},
\cite{cit_98500_Freud_book},
\cite{cit_20000_Suetin}, \cite{cit_3000_Chihara}, \cite{cit_5200_Nevai_book},~\cite{cit_5000_Ismail}. 
Let $\{ g_n (x) \}_{n=0}^\infty$, $\deg g_n = n$, be orthonormal polynomials on the real line, having positive leading coefficients:
\begin{equation}
\label{f1_22}
\int_\mathbb{R} g_n(x) g_m(x) d\mu_g(x) = \delta_{n,m},\qquad n,m\in\mathbb{Z}_+.
\end{equation}
Here $\mu_g$ is a (non-negative) measure on $\mathfrak{B}(\mathbb{R})$. It is well known that polynomials $g_n$ satisfy the following
recurrence relation:
\begin{equation}
\label{f1_25}
\widehat{a}_{n-1} g_{n-1}(x) + \widehat{b}_{n} g_{n}(x) +  \widehat{a}_{n} g_{n+1}(x) = x g_n(x),\qquad n\in \mathbb{Z}_+, 
\end{equation}
where $\widehat{a}_j>0$, $\widehat{b}_j\in\mathbb{R}$, $j\in\mathbb{Z}_+$; and $\widehat{a}_{-1}=0$, $g_{-1}(x)=0$.
Relation~(\ref{f1_25}) can be written in the following matrix form:
\begin{equation}
\label{f1_27}
G \vec g(x) = x \vec g(x),
\end{equation}
where $\vec g(x) = (g_0(x),g_1(x),...)^T$. Here $G$ is the Jacobi matrix with $\widehat{b}_k$s on the main diagonal, and with
$\widehat{a}_k$s on the first subdiagonal~(\cite{cit_2000_Akhiezer_book_1965}).
Polynomials
\begin{equation}
\label{f1_29}
K_n(t,x) = \sum_{k=0}^n g_k(t) g_k(x),\qquad n\in\mathbb{Z}_+,
\end{equation}
are called \textit{kernel polynomials} (see, e.g.,~\cite[p. 39]{cit_50000_Gabor_Szego} and observe that we use a slightly 
different notation). Kernel polynomials have \textit{the reproducing property} as well as other
important properties. In particular, if the measure support lies inside $(-\infty,b]$, and $t=t_0\geq b$, then kernel polynomials
$K_n(t_0,x)$ are orthogonal on $\mathbb{R}$ with respect to $(t_0-x)d\mu_g$, see~\cite[p. 40]{cit_50000_Gabor_Szego}.

Let $\vec c = \{ c_k \}_{k=0}^\infty$ be a row vector of arbitrary positive numbers. Consider the following polynomials:
\begin{equation}
\label{f1_31}
u_n(x) = u_n(\vec c; x) := \sum_{k=0}^n c_k g_k(x),\qquad n\in\mathbb{Z}_+.
\end{equation}
These polynomials generalize the above polynomials $K_n(t_0,x)$, since $c_k = g_k(t_0)>0$ in that case.
As we shall see, polynomials $u_n(x)$ need not to be orthogonal on the real line. However, they possess many interesting properties.
They will be called \textbf{modified kernel polynomials}.
Further examples of polynomials of type~(\ref{f1_31}) appear from Green's formula for solutions of the second-order
difference equation generated by~(\ref{f1_25})~\cite[p. 9]{cit_2000_Akhiezer_book_1965}. For example, 
one can choose $c_0=1$; $c_k = q_k(t_0)$, $k\geq 1$, when the latter numbers are positive.
Here $q_k$ are polynomials of the second kind for $g_k$. The use of Green's formula provides a compact representation for
these polynomials. Parameters $c_k$ may be also chosen to be some other solutions of the corresponding
second-order difference equation of $g_k$.
Other examples of modified kernel polynomials lead to some classical type Sobolev orthogonal polynomials.
The theory of Sobolev orthogonal polynomials is nowadays of a considerable interest.
A brief account of the theory state can be found in~\cite{cit_3500_CS_MB_JAT_2011} and a more recent survey is given in~\cite{cit_5150_M_X}.
Before we shall discuss the examples of Sobolev orthogonal polynomials related to~(\ref{f1_31}),
we  need the following definition.

\begin{definition} (\cite{cit_500460_Zagorodnyuk_UMZh_2017})
\label{d1_1}
A set $\Theta = \left(
J_3, J_5, \alpha, \beta
\right)$,
where $\alpha>0$, $\beta\in\mathbb{R}$, $J_3$ is a Jacobi matrix and
$J_5$ is a semi-infinite real symmetric five-diagonal matrix with positive numbers on the second subdiagonal,
is said to be
\textbf{a Jacobi-type pencil (of matrices)}.
\end{definition}
From this definition it follows that matrices $J_3$ and $J_5$ have the following form:
\begin{equation}
\label{ff1_5}
J_3 =
\left(
\begin{array}{cccccc}
b_0 & a_0 & 0 & 0 & 0 & \cdots\\
a_0 & b_1 & a_1 & 0 & 0 & \cdots\\
0 & a_1 & b_2 & a_2 & 0 &\cdots\\ 
\vdots & \vdots & \vdots & \ddots \end{array}
\right),\qquad a_k>0,\ b_k\in\mathbb{R},\ k\in\mathbb{Z}_+;
\end{equation}

\begin{equation}
\label{ff1_10}
J_5 =
\left(
\begin{array}{ccccccc}
\alpha_0 & \beta_0 & \gamma_0 & 0 & 0 & 0 & \cdots\\
\beta_0 & \alpha_1 & \beta_1 & \gamma_1 & 0 & 0 & \cdots\\
\gamma_0 & \beta_1 & \alpha_2 & \beta_2 & \gamma_2 & 0 & \cdots\\ 
0 & \gamma_1 & \beta_2 & \alpha_3 & \beta_3 & \gamma_3 & \cdots \\
\vdots & \vdots & \vdots &\vdots & \ddots \end{array}
\right),\ \alpha_n,\beta_n\in\mathbb{R},\ \gamma_n>0,\ n\in\mathbb{Z}_+.
\end{equation}

With a Jacobi-type pencil of matrices $\Theta$ one associates a system of polynomials
$\{ p_n(\lambda) \}_{n=0}^\infty$, such that
\begin{equation}
\label{ff1_15}
p_0(\lambda) = 1,\quad p_1(\lambda) = \alpha\lambda + \beta,
\end{equation}
and
\begin{equation}
\label{ff1_20}
(J_5 - \lambda J_3) \vec p(\lambda) = 0,
\end{equation}
where $\vec p(\lambda) = (p_0(\lambda), p_1(\lambda), p_2(\lambda),\cdots)^T$. 
Polynomials $\{ p_n(\lambda) \}_{n=0}^\infty$ are said to be \textit{associated to the Jacobi-type pencil of matrices $\Theta$}.

Relation~(\ref{ff1_20}) shows that $\vec p(\lambda)$ is a generalized eigenvector of a linear operator pencil 
(or, in other words, of an operator polynomial)
$J_5 - \lambda J_3$,
corresponding to eigenvalue $\lambda$.
Relation~(\ref{ff1_20}) is a generalization of relation~(\ref{f1_27}). This is a discretization of a $4$-th order differential
equation appearing in some physical problems~\cite{cit_500470_Zagorodnyuk_SIGMA_2017}. 
One can rewrite relation~(\ref{ff1_20}) in the following scalar form:
$$ \gamma_{n-2} p_{n-2}(\lambda) + (\beta_{n-1}-\lambda a_{n-1}) p_{n-1}(\lambda) + (\alpha_n-\lambda b_n) p_n(\lambda) +
$$
\begin{equation}
\label{ff1_30}
+ (\beta_n-\lambda a_n) p_{n+1}(\lambda) + \gamma_n p_{n+2}(\lambda) = 0,\qquad n\in\mathbb{Z}_+,
\end{equation}
where $p_{-2}(\lambda) = p_{-1}(\lambda) = 0$, $\gamma_{-2} = \gamma_{-1} = a_{-1} = \beta_{-1} = 0$.
Polynomials $p_n$ satisfy some special orthogonality relations as well as other interesting properties~\cite{cit_500480_Zagorodnyuk_JDEA_2018}.

Another possible generalization of~(\ref{f1_27}) of the following form:
\begin{equation}
\label{f1_48}
J_{2N+1} \vec p(\lambda) = \lambda^N \vec p(\lambda),
\end{equation}
with $J_{2N+1}$ being a $(2N+1)$-diagonal complex semi-infinite symmetric matrix,
appeared in the study of orthogonal polynomials on radial rays in the complex plane, 
see, e.g.,~\cite{cit_3700_Duran_Van_Assche_1996},\cite{cit_5170_Milovanovic},\cite{cit_3200_Choque_Rivero_Zagorodnyuk_2009} and
references therein for other generalizations.
Observe that pencils of threediagonal matrices and associated biorthogonal rational functions were studied in~\cite{cit_98000_Zhedanov_JAT}.

Let us briefly describe the content of the present paper.
The modified kernel polynomials satisfy a recurrence relation of type~(\ref{ff1_30}) (Theorem~\ref{t2_1}).
In order to prove this we emborder relation~(\ref{f1_27}) by two-diagonal matrices.
This idea for producing systems of polynomials satisfying~(\ref{ff1_20}) was proposed by one of the referees of~\cite{cit_500460_Zagorodnyuk_UMZh_2017}.
We also needed to change the sign of matrices to fit into~(\ref{ff1_20}) and to choose special two-diagonal matrices.

When $g_k$ are chosen to be Jacobi or Laguerre polynomials, these special cases have additional advantages.
Suitable choices of parameters $c_k$ imply $u_n$ to be Sobolev orthogonal polynomials.
Here we used an idea from~\cite{cit_500500_Zagorodnyuk_JAT_2020}.

A further selection of parameters $c_k$ gives some differential equations for $u_n$.
The above results will be gathered in Theorems~2 and~3.
Since $u_n$ are (generalized) eigenfunctions both for a difference operator and for a differential operator,
they may be viewed as classical type Sobolev orthogonal polynomials.

\noindent
{\bf Notations. }
As usual, we denote by $\mathbb{R}, \mathbb{C}, \mathbb{N}, \mathbb{Z}, \mathbb{Z}_+$,
the sets of real numbers, complex numbers, positive integers, integers and non-negative integers,
respectively. The superscript $T$ means the matrix transpose.
By $\mathbb{Z}_{k,l}$ we mean all integers $j$ satisfying the following inequality:
$k\leq j\leq l$; ($k,l\in\mathbb{Z}$).
By $\mathfrak{B}(\mathbb{R})$ we mean the set of all Borel subsets of $\mathbb{R}$.
By $\mathbb{P}$ we denote the set of all polynomials with complex coefficients.
For a complex number $c$ we denote
$(c)_0 = 1$, $(c)_1=c$, $(c)_k = c(c+1)...(c+k-1)$, $k\in\mathbb{N}$ (\textit{the shifted factorial or Pochhammer symbol}).
The generalized hypergeometric function is denoted by
$$ {}_m F_n(a_1,...,a_m; b_1,...,b_n;x) = \sum_{k=0}^\infty \frac{(a_1)_k ... (a_m)_k}{(b_1)_k ... (b_n)_k} \frac{x^k}{k!}, $$
where $m,n\in\mathbb{N}$, $a_j,b_l\in\mathbb{C}$.
By $\Gamma(z)$ and $\mathrm{B}(z)$ we denote the gamma function and the beta function, respectively.
By $J_\nu(z)$ we mean the Bessel function of the first kind.

\section{Modified kernel polynomials.}

Our first aim is to provide a recurrence relation for the modified kernel polynomials~(\ref{f1_31}).
\begin{theorem}
\label{t2_1}
Let $\{ g_n (x) \}_{n=0}^\infty$ ($\deg g_n = n$) be orthonormal polynomials on the real line with positive leading coefficients and
satisfying relation~(\ref{f1_25}).
Let $\{ c_k \}_{k=0}^\infty$ be a set of arbitrary positive numbers, and $u_n$ be the
modified kernel polynomials from~(\ref{f1_31}). 
The polynomials
\begin{equation}
\label{f2_2}
p_n(x) := \frac{1}{c_0 g_0} u_n(x) = \frac{1}{c_0 g_0} \sum_{j=0}^n c_j g_j(x),\qquad n\in\mathbb{Z}_+,
\end{equation}
are associated polynomials to the following Jacobi-type pencil:
\begin{equation}
\label{f2_4}
\widetilde\Theta = \left(
\mathbf{J}_3, \mathbf{J}_5, \widetilde\alpha, \widetilde\beta
\right), 
\end{equation}
where
\begin{equation}
\label{f2_6}
\widetilde\alpha = \frac{ c_1 }{ c_0 \widehat{a}_0 },\quad 
\widetilde\beta = 1 - \frac{ c_1 \widehat{b}_0 }{ c_0 \widehat{a}_0 },
\end{equation}
matrices $\mathbf{J}_3$, $\mathbf{J}_5$, have forms~(\ref{ff1_5}) and (\ref{ff1_10}), respectively, with
\begin{equation}
\label{f2_7}
a_n = \frac{1}{ c_{n+1}^2 },\quad b_n = -\frac{1}{ c_{n}^2 } - \frac{1}{ c_{n+1}^2 },
\end{equation}
$$ \alpha_n = 
\frac{ 2 \widehat{a}_{n} }{ c_n c_{n+1} } - \frac{ \widehat{b}_{n} }{ c_n^2 } -
\frac{ \widehat{b}_{n+1} }{ c_{n+1}^2 },\quad
\beta_n = 
\frac{ \widehat{b}_{n+1} }{ c_{n+1}^2 } - \frac{ \widehat{a}_{n+1} }{ c_{n+1} c_{n+2} } -
\frac{ \widehat{a}_{n} }{ c_n c_{n+1} }, $$
\begin{equation}
\label{f2_8}
\gamma_n = \frac{ \widehat{a}_{n+1} }{ c_{n+1} c_{n+2} },\qquad n\in\mathbb{Z}_+.
\end{equation}

Polynomials $p_n$ and $u_n$ satisfy relation~(\ref{ff1_30}) with the above coefficients from~(\ref{f2_8}).
\end{theorem}
\textbf{Proof.}
Let $\{ g_n (x) \}_{n=0}^\infty$, $\deg g_n = n$, be orthonormal polynomials on the real line with positive leading coefficients
and relations~(\ref{f1_22}),(\ref{f1_25}),(\ref{f1_27}) hold.
Define the modified kernel polynomials by relation~(\ref{f1_31}),
where $\{ c_k \}_{k=0}^\infty$ be a set of arbitrary positive numbers.
Consider the following semi-infinite two-diagonal matrix:
\begin{equation}
\label{f2_10}
C = 
\left(
\begin{array}{ccccc}
\frac{1}{c_0} & 0 & 0 & 0 & \cdots\\
-\frac{1}{c_1} & \frac{1}{c_1} & 0 & 0 & \cdots\\
0 & -\frac{1}{c_2} & \frac{1}{c_2} & 0 & \cdots\\
0 & 0 & -\frac{1}{c_3} & \frac{1}{c_3} & \cdots\\
\vdots & \vdots & \vdots & \vdots & \ddots\end{array}
\right).
\end{equation}
Observe that
\begin{equation}
\label{f2_15}
\vec g(x) = C \vec u(x), 
\end{equation}
where $\vec g(x) = (g_0(x),g_1(x),...)^T$, $\vec u(x) = (u_0(x),u_1(x),...)^T$.
Relation~(\ref{f2_15}) can be directly verified by relation~(\ref{f1_31}).
Substitute for $\vec g(x)$ from~(\ref{f2_15}) into~(\ref{f1_27}) to get
\begin{equation}
\label{f2_20}
(G C) \vec u(x) = x C \vec u(x).
\end{equation}
Denote by $C^*$ the formal adjoint of $C$:
\begin{equation}
\label{f2_25}
C^* = 
\left(
\begin{array}{ccccc}
\frac{1}{c_0} & -\frac{1}{c_1} & 0 & 0 & \cdots\\
0 & \frac{1}{c_1} & -\frac{1}{c_2} & 0 & \cdots\\
0 & 0 & \frac{1}{c_2} & -\frac{1}{c_3} & \cdots\\
\vdots & \vdots & \vdots & \vdots & \ddots\end{array}
\right).
\end{equation}
Then
\begin{equation}
\label{f2_30}
( - C^* G C ) \vec u(x) = x ( - C^* C ) \vec u(x).
\end{equation}
Denote
\begin{equation}
\label{f2_35}
\mathbf{J}_3 = - C^* C,\quad \mathbf{J}_5 =  - C^* G C.
\end{equation}
Let us calculate the elements of matrices $\mathbf{J}_3$ and $\mathbf{J}_5$.
Denote by $\vec e_n$ the semi-infinite column vector $(\delta_{j,n})_{j=0}^\infty$, $n\in\mathbb{Z}_+$.
Let $\vec e_{-1} =\vec e_{-2}$ be zero column vectors.
Observe that
\begin{equation}
\label{f2_37}
C \vec e_n = \frac{1}{c_n} \vec e_n - \frac{1}{ c_{n+1} } \vec e_{n+1}, \qquad n\in\mathbb{Z}_+;
\end{equation}
\begin{equation}
\label{f2_39}
C^* \vec e_m = -\frac{1}{c_m} \vec e_{m-1} + \frac{1}{ c_{m} } \vec e_{m}, \qquad m\in\mathbb{Z}_+.
\end{equation}
Then
\begin{equation}
\label{f2_41}
\mathbf{J}_3 \vec e_n = \frac{1}{c_n^2} \vec e_{n-1} -
\left(
\frac{1}{c_n^2} + \frac{1}{ c_{n+1}^2 }  
\right)
\vec e_{n} + \frac{1}{ c_{n+1}^2 } \vec e_{n+1}, \qquad n\in\mathbb{Z}_+.
\end{equation}
Notice that
\begin{equation}
\label{f2_43}
G \vec e_k = \widehat{a}_{k-1} \vec e_{k-1} + \widehat{b}_k \vec e_k + \widehat{a}_k \vec e_{k+1}, \qquad k\in\mathbb{Z}_+.
\end{equation}
By~(\ref{f2_37}),(\ref{f2_39}) and (\ref{f2_43}) it follows that
$$ \mathbf{J}_5 \vec e_n =
\frac{ \widehat{a}_{n-1} }{ c_n c_{n-1} } \vec e_{n-2} +
\left(
\frac{ \widehat{b}_{n} }{ c_n^2 } - \frac{ \widehat{a}_{n} }{ c_n c_{n+1} } -
\frac{ \widehat{a}_{n-1} }{ c_n c_{n-1} }
\right)
\vec e_{n-1} + $$
$$ +
\left(
\frac{ 2 \widehat{a}_{n} }{ c_n c_{n+1} } - \frac{ \widehat{b}_{n} }{ c_n^2 } -
\frac{ \widehat{b}_{n+1} }{ c_{n+1}^2 }
\right)
\vec e_n +
\left(
\frac{ \widehat{b}_{n+1} }{ c_{n+1}^2 } - \frac{ \widehat{a}_{n+1} }{ c_{n+1} c_{n+2} } -
\frac{ \widehat{a}_{n} }{ c_n c_{n+1} }
\right)
\vec e_{n+1} + $$
\begin{equation}
\label{f2_45}
+ \frac{ \widehat{a}_{n+1} }{ c_{n+1} c_{n+2} } \vec e_{n+2},\qquad n\in\mathbb{Z}_+,
\end{equation}
where $c_{-1} := 1$, $\widehat a_{-1} := 0$.
By~(\ref{f2_41}),(\ref{f2_45}) we see that matrices $\mathbf{J}_3, \mathbf{J}_5$ have forms~(\ref{ff1_5}), (\ref{ff1_10}), respectively.
Their entries are given by relations~(\ref{f2_7}) and (\ref{f2_8}).
Define polynomials $p_n$ by~(\ref{f2_2}). 
Then $p_0=1$, and $p_1$ can be directly calculated to get $\widetilde\alpha$ and $\widetilde\beta$ from~(\ref{f2_6}).
By~(\ref{f2_30}) we conclude that $p_n$ are associated to a Jacobi-type pencil $\left(
\mathbf{J}_3, \mathbf{J}_5, \widetilde\alpha, \widetilde\beta
\right)$.
$\Box$

Recall that the Jacobi polynomials $P_n^{(\alpha,\beta)}(x)$:
$$ P_n^{(\alpha,\beta)}(x) =
\left(
\begin{array}{cc} n+\alpha\\
n\end{array}
\right)
{}_2 F_1\left(-n,n+\alpha+\beta+1;\alpha+1;\frac{1-x}{2}\right),\  n\in\mathbb{Z}_+, $$
are orthogonal on $[-1,1]$ with respect to the weight $w(x) = (1-x)^\alpha (1+x)^\beta$,
$\alpha,\beta > -1$,~(\cite{cit_50_Bateman}). 
The orthonormal polynomials have the following form:
$$ \widehat P_0^{(\alpha,\beta)}(x) = \frac{1}{ \sqrt{ 2^{\alpha+\beta+1} \mathrm{B}(\alpha+1,\beta+1) } }, $$
$$ \widehat P_n^{(\alpha,\beta)}(x) = 
\sqrt{
\frac{ (2n+\alpha+\beta+1) n! \Gamma(n+\alpha+\beta+1) }{  2^{\alpha+\beta+1} \Gamma(n+\alpha+1) \Gamma(n+\beta+1)  }
}
P_n^{(\alpha,\beta)}(x),\quad n\in\mathbb{N}. $$
Polynomials $P_n^{(\alpha,\beta)}(x)$ and $\widehat P_n^{(\alpha,\beta)}(x)$ are solutions to the following
differential equation:
\begin{equation}
\label{f2_50}
(1-x^2) y'' + [ \beta-\alpha - (\alpha+\beta+2) x ] y' + n(n+\alpha+\beta+1) y = 0.
\end{equation}
Let $c>0$ be an arbitrary positive number.
Rewrite equation~(\ref{f2_50}) in the following form:
\begin{equation}
\label{f2_52}
D_{\alpha,\beta,c} y(x) = l_{n,c} y(x),\qquad n\in\mathbb{Z}_+,
\end{equation}
where
\begin{equation}
\label{f2_54}
D_{\alpha,\beta,c} := (x^2-1) \frac{d^2}{dx^2} + [ (\alpha+\beta+2) x + \alpha - \beta ] \frac{d}{dx} + c,
\end{equation}
\begin{equation}
\label{f2_56}
l_{n,c} := c + n(n+\alpha+\beta+1).
\end{equation}
It is important for us that $l_{n.c}$ are positive numbers.
On the other hand, relations~(\ref{f2_52})-(\ref{f2_56}) can be written for all real values of the parameter.
Define the following polynomials:
\begin{equation}
\label{f2_58}
P_n(\alpha,\beta,c,t_0;x) := \sum_{k=0}^n \frac{1}{ l_{k,c} } \widehat P_k^{(\alpha,\beta)}(t_0)
\widehat P_k^{(\alpha,\beta)}(x),\qquad n\in\mathbb{Z}_+,
\end{equation}
where $t_0\geq 1$ is an arbitrary parameter.
Since $\frac{1}{ l_{k,c} } \widehat P_k^{(\alpha,\beta)}(t_0) > 0$, we conclude that $P_n(\alpha,\beta,c,t_0;x)$ are
modified kernel polynomials of type~(\ref{f1_31}). Observe that scaled by eigenvalues polynomial kernels of
some \textit{Sobolev orthogonal polynomials} already appeared in the literature, see~\cite{cit_5100_LMMW_JAT_2018}.
Observe that
\begin{equation}
\label{f2_60}
D_{\alpha,\beta,c} P_n(\alpha,\beta,c,t_0;x) = \sum_{k=0}^n \widehat P_k^{(\alpha,\beta)}(t_0)
\widehat P_k^{(\alpha,\beta)}(x) =: p_n(\alpha,\beta,t_0;x),\quad n\in\mathbb{Z}_+.
\end{equation}
Polynomials $p_n(\alpha,\beta,t_0;x)$ are usual kernel polynomials. As it was noticed in the Introduction, they are
orthogonal on $\mathbb{R}$. Thus, polynomials $P_n(\alpha,\beta,c,t_0;x)$ fit into the scheme of paper~\cite{cit_500500_Zagorodnyuk_JAT_2020}
(see Condition~1 on page~3 therein). This means that $P_n(\alpha,\beta,c,t_0;x)$ are Sobolev orthogonal polynomials.

The case $t_0=1$ has a special advantage. In fact, $p_n(\alpha,\beta,1;x)$ is a multiple of the Jacobi polynomial for indices $(\alpha+1,\beta)$.
In this case, the differential equation for $p_n(\alpha,\beta,1;x)$, together with~(\ref{f2_60}), gives a differential
equation for $P_n(\alpha,\beta,c,1;x)$.

\begin{theorem}
\label{t2_2}
Let $\alpha,\beta>-1$; $c>0$, and $t_0\geq 1$, be arbitrary parameters.
Polynomials $P_n(x) = P_n(\alpha,\beta,c,t_0;x)$, given by relation~(\ref{f2_58}), are Sobolev orthogonal polynomials on $\mathbb{R}$:
$$ \int_{-1}^1 (P_n(x),P_n'(x),P_n''(x)) M_{\alpha,\beta,c}(x) \left(
\begin{array}{ccc} P_m(x) \\
 P_m'(x) \\
 P_m''(x) \end{array}
 \right)
(t_0-x) (1-x)^\alpha (1+x)^\beta dx  = $$
\begin{equation}
\label{f2_62}
= A_n \delta_{n.m},\qquad n,m\in\mathbb{Z}_+,
\end{equation}
where $A_n$ are some positive numbers and
$$ M_{\alpha,\beta,c} = $$
\begin{equation}
\label{f2_64}
=
\left(
\begin{array}{ccc} c \\
(\alpha+\beta+2) x + \alpha - \beta \\
x^2-1\end{array}
\right)
\left(
c, (\alpha+\beta+2) x + \alpha - \beta, x^2-1
\right).
\end{equation}
Moreover, the following differential equation holds for $P_n(\alpha,\beta,c,1;x)$:
\begin{equation}
\label{f2_67}
D_{\alpha+1,\beta,0} D_{\alpha,\beta,c} P_n(\alpha,\beta,c,1;x) = l_{n,0} D_{\alpha,\beta,c} P_n(\alpha,\beta,c,1;x),\quad n\in\mathbb{Z}_+,
\end{equation}
where $D_{\alpha,\beta,c}$, $l_{n,c}$ are defined as in~(\ref{f2_54}),(\ref{f2_56}).
\end{theorem}
\textbf{Proof.} All statements of the theorem readily follow from considerations before its formulation. 
$\Box$

\begin{example}
\label{e2_1}
Consider a polynomial of the following form:
\begin{equation}
\label{f2_70}
w_2(c;x) = \sum_{k=0}^2 \frac{ \widetilde{a}_k }{ \widetilde{b}_k + c } \widetilde g_k(x),\qquad c>0,
\end{equation}
where $\widetilde{a}_k > 0$, $\widetilde{b}_k \geq 0$, and
$\widetilde{g}_k$ is a real polynomial of degree~$k$ with a positive leading coefficient; $k=0,1,2$.
Of course, the polynomial $P_2(\alpha,\beta,c,t_0;x)$, with some fixed admissible parameters $\alpha,\beta,t_0$,
has form~(\ref{f2_70}). Another polynomial of form~(\ref{f2_70}) will appear below.
Let
\begin{equation}
\label{f2_72}
\widetilde{g}_k = \sum_{j=0}^k \mu_{k;j} x^j,\qquad k=0,1,2, 
\end{equation}
with $\mu_{k;j}\in\mathbb{R}$, $\mu_{k;k}>0$.
Substitute for $\widetilde{g}_k$ from~(\ref{f2_72}) into~(\ref{f2_70}) to get
$$ w_2(c;x) = \frac{ \widetilde{a}_2 }{ \widetilde{b}_2 + c } \mu_{2;2} x^2 +  \left(
\frac{ \widetilde{a}_1 }{ \widetilde{b}_1 + c } \mu_{1;1} + \frac{ \widetilde{a}_2 }{ \widetilde{b}_2 + c } \mu_{2;1}
\right) x + $$
\begin{equation}
\label{f2_74}
+ \frac{ \widetilde{a}_0 }{ \widetilde{b}_0 + c } \mu_{0;0} +
\frac{ \widetilde{a}_1 }{ \widetilde{b}_1 + c } \mu_{1;0} +
\frac{ \widetilde{a}_2 }{ \widetilde{b}_2 + c } \mu_{2;0}.
\end{equation}
The discriminant $D$ of the latter quadratic equation is equal to
$$ D = 
\left(
\frac{ \widetilde{a}_1 }{ \widetilde{b}_1 + c } \mu_{1;1} + \frac{ \widetilde{a}_2 }{ \widetilde{b}_2 + c } \mu_{2;1}
\right)^2 -
$$
\begin{equation}
\label{f2_76}
-
4
\frac{ \widetilde{a}_2 }{ \widetilde{b}_2 + c } \mu_{2;2}
\left(
\frac{ \widetilde{a}_0 }{ \widetilde{b}_0 + c } \mu_{0;0} +
\frac{ \widetilde{a}_1 }{ \widetilde{b}_1 + c } \mu_{1;0} +
\frac{ \widetilde{a}_2 }{ \widetilde{b}_2 + c } \mu_{2;0}
\right).
\end{equation}
In the case of $P_2(\alpha,\beta,c,t_0;x)$ we have $\widetilde{b}_0 = 0$, while $\widetilde{b}_1, \widetilde{b}_2 > 0$.
Therefore, for some small positive values of $c$, the discriminant $D$ is negative.
This fact shows that polynomials $P_n(\alpha,\beta,c,t_0;x)$ (for some small values of $c$) are not orthogonal polynomials on the real line.

\end{example}

Recall that the (generalized) Laguerre polynomials $L_n^{\alpha}(x)$:
$$ L_n^{\alpha}(x) =
\left(
\begin{array}{cc} n+\alpha\\
n\end{array}
\right)
{}_1 F_1\left( -n;\alpha+1;x \right),\qquad  n\in\mathbb{Z}_+, $$
are orthogonal on $[0,+\infty)$ with respect to the weight $w(x) = x^\alpha e^{-x}$,
$\alpha > -1$,~(\cite{cit_50_Bateman}). 
Orthonormal polynomials $\widehat L_n^{\alpha}(x)$ have the following form:
$$ \widehat L_n^{\alpha}(x) = \frac{ (-1)^n }{ \sqrt{ \Gamma(\alpha+1) \left(
\begin{array}{cc} n+\alpha\\
n\end{array}
\right)  } } L_n^{\alpha}(x),\qquad  n\in\mathbb{Z}_+. $$
Polynomials $L_n^{\alpha}(x)$ satisfy the following differential equation:
\begin{equation}
\label{f2_78}
x y'' + ( \alpha + 1 - x ) y' = -ny.
\end{equation}
Polynomials
$$ \mathcal{L}_n(x) = \mathcal{L}_n(\alpha;x) :=  \frac{ 1 }{ \sqrt{ \Gamma(\alpha+1) \left(
\begin{array}{cc} n+\alpha\\
n\end{array}
\right)  } } L_n^{\alpha}(-x),\qquad  n\in\mathbb{Z}_+, $$
are orthonormal on $(-\infty,0]$ with respect to the weight $(-x)^\alpha e^x$, see~\cite[section 2.3]{cit_50000_Gabor_Szego}.
By~(\ref{f2_78}) it follows that $\mathcal{L}_n(x)$ are solutions of the following differential equation:
\begin{equation}
\label{f2_80}
x \mathcal{L}_n''(x) + (\alpha+1+x) \mathcal{L}_n'(x) = n \mathcal{L}_n(x).
\end{equation}

Let $c>0$ be an arbitrary positive number.
Equation~(\ref{f2_80}) can be written in the following form:
\begin{equation}
\label{f2_82}
D_{\alpha,c} \mathcal{L}_n(x) = \lambda_{n,c} \mathcal{L}_n(x),\qquad n\in\mathbb{Z}_+,
\end{equation}
where
\begin{equation}
\label{f2_84}
D_{\alpha,c} := x \frac{d^2}{dx^2} + (\alpha+1+x) \frac{d}{dx} + c,
\end{equation}
\begin{equation}
\label{f2_86}
\lambda_{n,c} := c + n.
\end{equation}
Observe that $\lambda_{n,c}$ are positive numbers.
Define the following polynomials:
\begin{equation}
\label{f2_88}
L_n(\alpha,c,t_0;x) := \sum_{k=0}^n \frac{1}{ \lambda_{k,c} } \mathcal{L}_k(\alpha;t_0) \mathcal{L}_k(\alpha;x),\qquad n\in\mathbb{Z}_+,
\end{equation}
where $t_0\geq 0$ is an arbitrary parameter.
Since $\frac{1}{ \lambda_{k,c} } \mathcal{L}_k(\alpha;t_0) > 0$, then $L_n(\alpha,c,t_0;x)$ are
modified kernel polynomials of type~(\ref{f1_31}). 
Notice that
\begin{equation}
\label{f2_90}
D_{\alpha,c} L_n(\alpha,c,t_0;x) = \sum_{k=0}^n \mathcal{L}_k(\alpha;t_0) \mathcal{L}_k(\alpha;x) =: p_n(\alpha,t_0;x),\quad n\in\mathbb{Z}_+.
\end{equation}
Polynomials $p_n(\alpha,t_0;x)$ are usual kernel polynomials. Since $t_0\geq 0$, they are
orthogonal on the real line. Thus, polynomials $L_n(\alpha,c,t_0;x)$ fit into the scheme of paper~\cite{cit_500500_Zagorodnyuk_JAT_2020} and
therefore $L_n(\alpha,c,t_0;x)$ are Sobolev orthogonal polynomials.

The case $t_0=0$ has a special interest. In fact, $p_n(\alpha,0;x)$ is a constant multiple of $\mathcal{L}_n(\alpha+1;x)$.
In this case, the differential equation for $p_n(\alpha,0;x)$, together with~(\ref{f2_80}), gives a differential
equation for $L_n(\alpha,c,0;x)$. 

\begin{theorem}
\label{t2_3}
Let $\alpha>-1$; $c>0$, and $t_0\geq 0$, be arbitrary parameters.
Polynomials $L_n(x) = L_n(\alpha,c,t_0;x)$, given by relation~(\ref{f2_88}), are Sobolev orthogonal polynomials on $\mathbb{R}$:
$$ \int_{-\infty}^0 (L_n(x),L_n'(x),L_n''(x)) M_{\alpha,c}(x) \left(
\begin{array}{ccc} L_m(x) \\
 L_m'(x) \\
 L_m''(x) \end{array}
 \right)
(t_0-x) (-x)^\alpha e^{x} dx  = $$
\begin{equation}
\label{f2_92}
= B_n \delta_{n.m},\qquad n,m\in\mathbb{Z}_+,
\end{equation}
where $B_n$ are some positive numbers and
\begin{equation}
\label{f2_94}
M_{\alpha,c} =
\left(
\begin{array}{ccc} c \\
\alpha + 1 + x \\
x\end{array}
\right)
\left(
c, \alpha + 1 + x, x
\right).
\end{equation}
Moreover, the following differential equation holds for $L_n(\alpha,c,0;x)$:
\begin{equation}
\label{f2_97}
D_{\alpha+1,0} D_{\alpha,c} L_n(\alpha,c,0;x) = \lambda_{n,0} D_{\alpha,c} L_n(\alpha,c,0;x),\quad n\in\mathbb{Z}_+,
\end{equation}
where $D_{\alpha,c}$, $\lambda_{n,c}$ are defined as in~(\ref{f2_84}),(\ref{f2_86}).
For an arbitrary $x<0$, and $c\in\mathbb{N}$, polynomials $L_n(\alpha,c,0;x)$ admit the following integral representation:
$$ L_n(\alpha,c,0;x) = \frac{1}{ \Gamma(\alpha+1) n! } e^{-x} (-x)^{ -\frac{\alpha}{2} } * $$
\begin{equation}
\label{f2_105}
* \int_0^\infty \int_0^t e^{-t} t^{ \frac{\alpha}{2} - c } \theta^{n+c-1} J_\alpha \left( 2 \sqrt{-tx} \right)
{}_2 F_0\left(
-n,1;-;-\frac{1}{\theta}
\right)
d\theta dt,\quad n\in\mathbb{Z}_+.
\end{equation}

\end{theorem}
\textbf{Proof.} All statements of the theorem, except the last one, follow from considerations before its formulation. 
Fix arbitrary $\alpha>-1$, $c\in\mathbb{N}$.
By the definitions of $L_n(\alpha,c,0;x), \mathcal{L}_k(\alpha;t_0)$, using $L_n^\alpha(0) = \frac{(\alpha+1)_n}{n!}$, we get
\begin{equation}
\label{f2_107}
L_n(\alpha,c,0;x) = \frac{1}{ \Gamma(\alpha+1) } \sum_{k=0}^n \frac{1}{k+c} L_k^\alpha(-x).
\end{equation}
If $x<0$, then~(see~\cite[p. 206]{cit_20000_Suetin})
\begin{equation}
\label{f2_109}
L_n^\alpha(-x) =
\frac{1}{ n! } e^{-x} (-x)^{ -\frac{\alpha}{2} }
\int_0^\infty e^{-t} t^{ n+\frac{\alpha}{2} } J_\alpha \left( 2 \sqrt{-tx} \right)
dt.
\end{equation}
Therefore
$$ L_n(\alpha,c,0;x) = \frac{1}{ \Gamma(\alpha+1) } 
e^{-x} (-x)^{ -\frac{\alpha}{2} } * $$
\begin{equation}
\label{f2_109_1}
*
\int_0^\infty
\left(
\sum_{k=0}^n \frac{1}{k+c} 
\frac{t^k}{ k! } 
\right)
e^{-t} t^{ \frac{\alpha}{2} } J_\alpha \left( 2 \sqrt{-tx} \right)
dt,\qquad x<0;\ n\in\mathbb{Z}_+.
\end{equation}

Consider the following polynomials of $t$:
$$ f_n(t) = f_n(c;t) := \sum_{k=0}^n \frac{1}{k+c} \frac{t^k}{k!},\qquad n\in\mathbb{Z}_+,\ c\in\mathbb{N}. $$
Observe that
\begin{equation}
\label{f2_110}
\left( f_n(t) t^c \right)' = t^{c-1} \widetilde f_n(t),\qquad t\in\mathbb{R}, 
\end{equation}
where 
$$ \widetilde f_n(t) := \sum_{k=0}^n \frac{t^k}{k!},\qquad n\in\mathbb{Z}_+. $$
The polynomial $\widetilde f_n(t)$ admit the following representation:
\begin{equation}
\label{f2_115}
\widetilde f_n(t) = \frac{1}{n!} t^n {}_2 F_0\left(
-n,1;-;-\frac{1}{t}
\right),
\end{equation}
where for $t=0$ we mean the limit value. Relation~(\ref{f2_115}) can be checked by writing ${}_2 F_0$ as a terminating series
and changing the index of summation.
Integrating~(\ref{f2_110}) and using~(\ref{f2_115}) we see that
\begin{equation}
\label{f2_117}
f_n(t) = \frac{1}{n!} t^{-c} \int_0^t \theta^{n+c-1} {}_2 F_0\left(
-n,1;-;-\frac{1}{\theta}
\right) d\theta,\qquad t>0.
\end{equation}
By relations~(\ref{f2_117}) and~(\ref{f2_109_1}) we obtain the required integral
representation~(\ref{f2_105}).
Notice that for $t=0$ formula~(\ref{f2_117}) is not valid. However, the value of the integral in~(\ref{f2_105}) does not
depend on this value. 
$\Box$

\begin{example}
\label{e2_2}
Observe that
$$ L_2(\alpha,c,t_0;x) := \sum_{k=0}^2 \frac{ \mathcal{L}_k(\alpha;t_0) }{ k+c } \mathcal{L}_k(\alpha;x), $$
has form~(\ref{f2_70}). We have $\widetilde{b}_0 = 0$, while $\widetilde{b}_1, \widetilde{b}_2 > 0$.
Thus, for some small positive values of $c$, the discriminant $D$ from~(\ref{f2_76}) is negative.
This fact shows that polynomials $L_n(\alpha,c,t_0;x)$ (at least for some values of $c$) are not orthogonal polynomials on the real line.

\end{example}

Consider polynomials $P_n(\alpha,\beta,c,1;x)$ (see~(\ref{f2_58})) with $\alpha=\beta=-\frac{1}{2}$. Denote by
$$ T_n(x) = \cos(n\arccos x),\qquad n\in\mathbb{Z}_+,\ x\in [-1,1], $$
the Chebyshev polynomials of the first kind. Orthonormal polynomials are given by (\cite{cit_20000_Suetin})
\begin{equation}
\label{f2_118}
\widehat T_0(x) = \frac{1}{ \sqrt{\pi} },\quad \widehat T_n(x) = \sqrt{ \frac{2}{\pi} } \cos(n\arccos x),\qquad n\in\mathbb{N},\ 
x\in [-1,1]. 
\end{equation}
Then
$$ t_n(c;x) = t_n(x) := P_n\left(-\frac{1}{2},-\frac{1}{2},c,1;x\right) = $$
\begin{equation}
\label{f2_119}
= \frac{1}{\pi c} + \frac{2}{\pi} \sum_{k=1}^n \frac{1}{k^2+c} T_k(x),\qquad n\in\mathbb{Z}_+;\ c>0;\ x\in\mathbb{C}.
\end{equation}
Here we used relations~(\ref{f2_58}),(\ref{f2_118}), and the equality $T_k(1)=1$.
The explicit representation~(\ref{f2_119}) shows that  
\begin{equation}
\label{f2_125}
|t_n(c;x)| \leq \frac{1}{\pi c} + \frac{2}{\pi} \sum_{k=1}^n \frac{1}{k^2}
\leq \frac{1}{\pi c} + \frac{2}{\pi}n,\qquad n\in\mathbb{Z}_+;\ c>0;\ x\in[-1,1].
\end{equation}
Since (\cite[p. 91]{cit_20000_Suetin})
$$ |T_n'(x)| \leq n^2,\qquad x\in[-1,1], $$
then
\begin{equation}
\label{f2_127}
|t_n'(c;x)| \leq \frac{2}{\pi} \sum_{k=1}^n \frac{k^2}{c + k^2}
\leq \frac{2}{\pi} n,\qquad n\in\mathbb{N};\ c>0;\ x\in[-1,1].
\end{equation}
Observe that
$$ \pi t_1(c;x) = \frac{1}{c} + \frac{2}{c+1} x, $$
$$ \pi t_2(c;x) = \frac{4}{c+4} x^2 + \frac{2}{c+1} x + \frac{1}{c} - \frac{2}{c+4}. $$
The root of $t_1$ is $x_1 = -\frac{c+1}{2c} \rightarrow -\frac{1}{2}$, as $c\rightarrow+\infty$.
We see that it need not to lie inside $[-1,1]$ for small values of $c$.
The roots of $t_2$ can be non-real for small $c$.
Further properties of polynomials $t_n(c;x)$ will be studied elsewhere.

\begin{center}
{\large\bf 
On modified kernel polynomials and classical type Sobolev orthogonal polynomials}
\end{center}
\begin{center}
{\bf S.M. Zagorodnyuk}
\end{center}

In this paper we study modified kernel polynomials:
$u_n(x) = \sum_{k=0}^n c_k g_k(x)$, depending on parameters $c_k>0$, where $\{ g_k \}_0^\infty$ are orthonormal
polynomials on the real line.
Besides kernel polynomials with $c_k = g_k(t_0)>0$, for example, $c_k$ may be chosen to be some other solutions of the corresponding
second-order difference equation of $g_k$.
It is shown that all these polynomials satisfy a $4$-th order recurrence relation.
The cases with $g_k$ being Jacobi or Laguerre polynomials are of a special interest.
Suitable choices of parameters $c_k$ imply $u_n$ to be Sobolev orthogonal polynomials with a $(3\times 3)$ matrix measure.
Moreover, a further selection of parameters gives differential equations for $u_n$.
In the latter case, polynomials $u_n(x)$ are solutions to a generalized eigenvalue problems both in $x$ and in $n$.

\vspace{1.5cm}

V. N. Karazin Kharkiv National University \newline\indent
School of Mathematics and Computer Sciences \newline\indent
Department of Higher Mathematics and Informatics \newline\indent
Svobody Square 4, 61022, Kharkiv, Ukraine

Sergey.M.Zagorodnyuk@gmail.com; Sergey.M.Zagorodnyuk@univer.kharkov.ua

}
\end{document}